# $6^m$ Theorem for Prime numbers

**Gandarawatta R.W.M.P.I.S.B.[1], Perera S.P.C.[2], Rathnayake R.M.L.S.[2]**

[1]Department of Engineering Mathematics, Faculty of Engineering, University of Peradeniya, Sri Lanka

[2]Department of Engineering Mathematics, Faculty of Engineering, University of Peradeniya, Sri Lanka

[2]Department of Electrical and Electronic Engineering, Faculty of Engineering, University of Peradeniya, Sri Lanka

issampath@yahoo.com

## Abstract

We show that for any $\boldsymbol{P = 6^{m+1} \bullet N - 1}$ is a **prime** number for any $\boldsymbol{1 < N \leq 13}$, $\boldsymbol{N \neq 8}$ and $\boldsymbol{N \neq i^{m+1}\ Mod(6i + 1)}$ where $i \in \boldsymbol{Z^+}$ and $m \in odd\ \boldsymbol{Z^+}$ for $\boldsymbol{1 < N \leq 13\ and\ N \neq 8}$ and also we further discuss that
$\boldsymbol{P = 6^{m+1} \bullet N - 1}$ is a **prime** number for $\boldsymbol{N} > 13$ if and only if,
$\boldsymbol{N \neq i^{m+1}\ Mod(6i + 1) + (6i + 1)a}$ ; $i, a \in \boldsymbol{Z^+}$

## 1. Introduction

In 300 B.C., Euclid proved that there are infinitely many prime numbers in the number line and till today, mathematicians all around the world are searching for different methods to find those. This quest has two main objectives, namely the first one is to find large prime numbers and some algorithms such as 'Mersenne Primes[1][2], Fermat Search[3]' are used for that. Secondly to extract consecutive primes, and 'Sieve of Eratosthenes[4]' is a commonly used algorithm for that purpose.

By 2018, Fermat numbers have been factorized from $F_o$ to $F_{11}$ and in 2016[5], Boklan and Conway proved that the existence of another Fermat prime number is less than one in a billion[6]. The largest prime number found so far is M77232917[7], which is discovered on 26.12.2017 as the 50th Mersenne prime number and it consists of more than 23 million digits. The searching mechanism used for the above purpose was Great Internet Mersenne Prime Search (GIMPS) which is introduced in 1996, but from 06.09.2008 to present day[6], it could find only 5 more prime numbers and we can see the probability of finding a Mersenne prime number is declined over time. Therefore, looking for new methods to find prime numbers is undoubtedly beneficial and we have studied a new subset of prime numbers as a result and propose the $6^m$ theorem with this.

The main advantage of the $6^m$ theorem is that it can directly show $\boldsymbol{P = 6^{m+1} \bullet N - 1}$ is a prime number when $\boldsymbol{N \neq i^{m+1}\ Mod(6i + 1)}$and $\boldsymbol{1 < N \leq 13}$ for any $m \in odd\ \boldsymbol{Z^+}$and $i \in \boldsymbol{Z^+}$. Also when $N > 13$, $\boldsymbol{P}$ can be identified as a prime or a composite number in around N/6 -2 calculations and it's a significant improvement over the current algorithms. We performed a case study for m=1 and proved that the answers for $i^2 Mod(6i + 1)$ is given from 6 equations only.

One of the main algorithms to find consecutive primes is 'Sieve of Eratosthenes' method and by using primary equations of $6^m$ theorem, we can show that our algorithm has lesser number of calculations for an any given number range than Sieve of Eratosthenes method. In both algorithms, the prime numbers are obtained by removing none prime numbers.

In Sieve of Eratosthenes algorithm, multiplications of 2 and 3 are removed at first, and the same principle is followed by the primary equations in 6m theorem by considering only 6n+1 and 6n-1 numbers for all cases and that saves a lot of calculations.

When removing the multiplications of any number to find prime numbers, removing the multiplications of that number from its' square value is sufficient. This can save a lot of calculations when findings prime in a large range. For example, when removing dividends of 10223 from a certain range, the first number to be removed is 102232=104,509,729 as 10223x2 is same as 2x10223 and has been removed in an earlier calculation. This phenomenon can be obtained from the primary equations of 6m theorem by considering i=j, j+1, j+2,... and so on.

Hence the primary equations of $6^m$ theorem is removing the n values of 6n+1 and 6n-1 series, the prime numbers can be found from this method are of 6 times the amount coming from the Sieve of Eratosthenes algorithm for a particular range.

## $6^m$ Theorem for Prime numbers

If $\boldsymbol{N} \neq \boldsymbol{i}^{\boldsymbol{m}+\boldsymbol{1}}\ \boldsymbol{Mod}(\boldsymbol{6i}+\boldsymbol{1})$ ; $i \in \boldsymbol{Z}^+$ $m \in odd\ \boldsymbol{Z}^+$ then,

For $\boldsymbol{1} < \boldsymbol{N} \leq \boldsymbol{13}\ \boldsymbol{and}\ \ \boldsymbol{N} \neq \boldsymbol{8}$, $\quad \boldsymbol{P} = \boldsymbol{6}^{\boldsymbol{m}+\boldsymbol{1}} \bullet \boldsymbol{N} - \boldsymbol{1}$ is a **Prime Number**.
For $\boldsymbol{N} > 13$, $\quad \boldsymbol{P} = \boldsymbol{6}^{\boldsymbol{m}+\boldsymbol{1}} \bullet \boldsymbol{N} - \boldsymbol{1}$ is a **Prime Number if and only if**
$\boldsymbol{N} \neq \boldsymbol{i}^{\boldsymbol{m}+\boldsymbol{1}}\ \boldsymbol{Mod}(\boldsymbol{6i}+\boldsymbol{1}) + (\boldsymbol{6i}+\boldsymbol{1})\boldsymbol{a}$ ; $i, a \in \boldsymbol{Z}^+$

## 2. Proof

### 2.1 Definition

Any positive integer $\boldsymbol{Z}^+$ can be represented as $\boldsymbol{Z}^+ = \{6n, 6n+1, 6n+2, 6n+3, 6n+4, 6n+5;$
$n = 0,1,2,3, \ldots\}$
Then any prime number $\mathbf{P}$ which is larger than 3 can be identified as; $\mathbf{P} \subset \{6n+1, 6n+5\ ;$
$n = 0,1,2,3, \ldots\}$
By simplifying further, $\mathbf{P} \subset \{6n-1, 6n+1\ ; n = 1,2,3, \ldots\}$
Two series $\mathbf{S1}$ and $\mathbf{S2}$ can be defined as; $\mathbf{S1} = \{6n-1\ ;\ n = 1,2,3, \ldots\}$, $\mathbf{S2} = \{6n+1\ ;$
$n = 1,2,3, \ldots\}$

**Case 1:** When $(6a-1) \in \mathbf{S1}, (6b-1) \in \mathbf{S1}$ then let's take $\boldsymbol{C1} = (6a-1)(6b-1)$
$\boldsymbol{C1} = (6a-1)(6b-1) = 6[6ab-a-b]+1 = 6n+1\ ;$ $n = 6ab-a-b > 0$ ; hence $\boldsymbol{C1} \in \mathbf{S2}$

$$a = \frac{n+b}{6b-1}$$

For $\boldsymbol{C1} \in \mathbf{P}$ then $a \notin \boldsymbol{Z}^+$ and then, $n+b \neq k(6b-1)\ ;\ k \in \boldsymbol{Z}^+$

$$n \neq -b + k(6b-1)$$

$$k = i, b = j \Rightarrow n \neq -j + i(6j-1)$$

$$n \neq 6ij - i - j \ \ldots\ldots\ldots\ldots\ldots (1) \quad \text{where } i, j \in \boldsymbol{Z}^+$$

**Case 2:** When $(6a+1) \in \mathbf{S2}, (6b+1) \in \mathbf{S2}$ then let's take $\boldsymbol{C2} = (6a+1)(6b+1)$
$\boldsymbol{C2} = (6a+1)(6b+1) = 6[6ab+a+b]+1 = 6n+1;$ $n = 6ab+a+b > 0$ ; hence $\boldsymbol{C2} \in \mathbf{S2}$

$$a = \frac{n-b}{6b+1}$$

For $\boldsymbol{C2} \in \mathbf{P}$ then $a \notin \boldsymbol{Z}^+$ and then, $n-b \neq k(6b+1)$ ; $k \in \boldsymbol{Z}^+$

$$n \neq b + k(6b+1)$$

$$k = i, b = j \Rightarrow n \neq j + i(6j+1)$$

$$n \neq 6ij + i + j \ \ldots\ldots\ldots\ldots\ldots (2) \quad \text{where } i, j \in \boldsymbol{Z}^+$$

**Case 3:**
When $(6a+1) \in \mathbf{S2}, (6b-1) \in \mathbf{S1}$ then let's take $\boldsymbol{C3} = (6a+1)(6b-1)$
$\boldsymbol{C3} = (6a+1)(6b-1) = 6[6ab-a+b]-1 = 6n-1\ ;$ $n = 6ab-a+b > 0$ ; hence $\boldsymbol{C3} \in \mathbf{S1}$

$$a = \frac{n-b}{6b-1}$$

For $\boldsymbol{C3} \in \mathbf{P}$ then $a \notin \boldsymbol{Z}^+$and then, $n - b \neq k(6b-1)$ ; $k \in \boldsymbol{Z}^+$

$$n \neq b + k(6b-1)$$

$$k = i, b = j \Rightarrow n \neq j + i(6j-1)$$

$$n \neq 6ij - i + j \text{ ...............③} \quad \text{where } i, j \in \boldsymbol{Z}^+$$

Similarly ⇒ $b = \frac{n+a}{6a+1}$

For $\boldsymbol{C3} \in \mathbf{P}$ then $b \notin \boldsymbol{Z}^+$ and then, $n + a \neq k(6a+1)$ ; $k \in \boldsymbol{Z}^+$

$$n \neq k(6a+1) - a$$

$$k = i, a = j \Rightarrow n \neq i(6j+1) - j$$

$$n \neq 6ij + i - j \text{ ...............④} \quad \text{where } i, j \in \boldsymbol{Z}^+$$

Same values for $n$ can be obtained from equations ③ and ④ by changing $i$ → $j$ and $j$ → $i$ , hence equation ③ is used for further analysis.

## 2.2 Analysis of S1

The following study is performed for $\mathbf{S1} = \{6n - 1\ ;\ n = 1,2,3, \ldots\}$ series only

From equation ③ , $n \neq 6ij - i + j$ ; $i, j \in \boldsymbol{Z}^+$

For $n = 6k$ ; $k \in \boldsymbol{Z}^+$ then $-i + j = 6a_1$ ; $a_1 \in \boldsymbol{Z_1} \subset \boldsymbol{Z}$ is a must

Then,

$$j = 6a_1 + i$$

$$n = 6i(6a_1 + i) + 6a_1$$

$$n = 6^2 a_1 i + 6(i^2 + a_1) = 6[6a_1 i + i^2 + a_1]$$

Even more if $n = 6^2 k$ , $k \in \boldsymbol{Z}^+$ ; $\quad i^2 + a_1 = 6a_2;\ a_2 \in \boldsymbol{Z_1} \subset \boldsymbol{Z}$

$$a_1 = 6a_2 - i^2$$

$$n = 6^2 i(6a_2 - i^2) + 6 \times 6a_2$$

$$n = 6^3 a_2 i - 6^2 i^3 + 6^2 a_2$$

$$n = 6^3 a_2 i + 6^2(-i^3 + a_2)$$

$$n = 6^2[6a_2 i + (-i^3 + a_2)]$$

Also if $n = 6^3 k, k \in \boldsymbol{Z}^+$; $\quad -i^3 + a_2 = 6a_3$

$$a_2 = 6a_3 + i^3; a_3 \in \boldsymbol{Z_1} \subset \boldsymbol{Z}$$

$$n = 6^3 i(6a_3 + i^3) + 6^2 \times 6a_3$$

$$n = 6^4 a_3 i + 6^3(i^4 + a_3)$$

$$n = 6^3[6a_3 i + i^4 + a_3]$$

Considering the pattern

For odd m ; $n = 6^m[6a_m i + i^{m+1} + a_m]$; $i \in \boldsymbol{Z}^+$, $a_m \in \boldsymbol{Z_1} \subset \boldsymbol{Z}$

$$n = 6^m[(6i+1)a_m + i^{m+1}] \text{ ...............⑤}$$

For even m; $n = 6^m[6a_m i - i^{m+1} + a_m]$; $i \in \boldsymbol{Z}^+, a_m \in \boldsymbol{Z_1} \subset \boldsymbol{Z}$

$$n = 6^m[(6i+1)a_m - i^{m+1}]; \text{ ,...............⑥}$$

For non primes in $\mathbf{S1}_{6^m} = \{6n - 1\ ;\ n = c \bullet 6^m\ ;\ \mathrm{m} \in \boldsymbol{Z}^+, \mathrm{c} \in \boldsymbol{Z}^+{}_1 \subset \boldsymbol{Z}^+\ \}$ ; $\mathbf{S1}_{6^m} \subset \mathbf{S1}$

$$\mathbf{S1}_{6^m} = \{\mathbf{S1}_{6^{m,odd}}, \mathbf{S1}_{6^{m,even}}\}$$

### 2.2.1 Analysis of $\mathbf{S1}_{6^{m,odd}}$

Considering only $\mathbf{S1}_{6^{m,odd}}$ series,

For odd m ; $n = 6^m[6a_m i + i^{m+1} + a_m]$

$$n = 6^m[(6i+1)a_m + i^{m+1}]\ ; i \in \boldsymbol{Z}^+, a_m \in \boldsymbol{Z} \text{...............⑤}$$

$n > 0$

From ⑤ $\quad 6^m[(6i+1)a_m + i^{m+1}] > 0$

$$(6i+1)a_m + i^{m+1} > 0$$

$$a_m > \frac{-i^{m+1}}{6i+1}$$

Therefore, $a_m = int_{CL}\left[\frac{-i^{m+1}}{6i+1}\right] + a$ ; $a \in \mathbf{0}, \boldsymbol{Z}^+$

$$a_m = -int_{FL}\left[\frac{i^{m+1}}{6i+1}\right] + a$$

Substituting back to ⑤; $n = 6^m\left[(6i+1)\left(-int_{FL}\left[\frac{i^{m+1}}{6i+1}\right] + a\right) + i^{m+1}\right]$

$$\boldsymbol{n} = 6^m\left[i^{m+1} - (6i+1)int_{FL}\left[\frac{i^{m+1}}{6i+1}\right] + (6i+1)a\right]$$

$i^{m+1} - (6i+1)int_{FL}\left[\frac{i^{m+1}}{6i+1}\right] = i^{m+1}\, Mod(6i+1)$ ............... fact

Therefore, any **non-prime** number in $\mathbf{S1}_{6^{\mathbf{m}},\mathbf{odd}}$ series can be found by,
$n = 6^m[i^{m+1}Mod(6i+1) + (6i+1)a]$ ; $i \in \boldsymbol{Z}^+$; $a \in \mathbf{0}, \boldsymbol{Z}^+$; $m \in odd\ \boldsymbol{Z}^+$

Hence, to generate any prime number **P** in $\mathbf{S1}_{6^{\mathbf{m}},\mathbf{odd}}$ series, $n$ must satisfy,
$n \neq 6^m[i^{m+1}Mod(6i+1) + (6i+1)a]$ ; $i \in \boldsymbol{Z}^+$; $a \in \mathbf{0}, \boldsymbol{Z}^+$; $m \in odd\ \boldsymbol{Z}^+$
Hence, $n \neq 6^m[i^{m+1}Mod(6i+1) + (6i+1)a] \neq 6^m$•N then,

$$\boldsymbol{P = 6 \bullet 6^m \bullet N - 1 = 6^{m+1} \bullet N - 1}$$

Any **non-prime** number in $\mathbf{S1}_{6^{\mathbf{m}},\mathbf{odd}}$ is given by,
$i^{m+1}Mod(6i+1) + (6i+1)a = \mathrm{N}$
when $a = 0$, $N = i^{m+1}Mod(6i+1)$
When $i = 1$, $i^{m+1}Mod(6i+1) = 1$ for any $m \in \boldsymbol{Z}^+$,
Therefore, at $i = 1$, $i^{m+1}Mod(6i+1) + (6i+1)a = 1 + 7a = \{1,8,15,22, \ldots.\} = \mathrm{N}$ ; $a \in \mathbf{0}, \boldsymbol{Z}^+$
and for $i \geq 2, a \geq 1$ , $N = i^{m+1}Mod(6i+1) + (6i+1)a > 13$

Hence, any Prime number in $\mathbf{S1}_{6^{\mathbf{m}},\mathbf{odd}}$ must satisfy,
$i^{m+1}Mod(6i+1) + (6i+1)a \neq \mathrm{N}$
For $a = 0$, $N \neq i^{m+1}Mod(6i+1)$
For $i = 1$, $\mathrm{N} \neq \{1,8,15,22, \ldots.\}$
For $i \geq 2, a \geq 1$ , $N \neq i^{m+1}Mod(6i+1) + (6i+1)a > 13$

Considering all the conditions above,
If $\mathrm{N} \neq i^{m+1}Mod(6i+1)$; $i \in \boldsymbol{Z}^+, m \in odd\ \boldsymbol{Z}^+$
Then, for $1 < N \leq 13$ $and$ $N \neq 8$, $\boldsymbol{P = 6^{m+1} \bullet N - 1}$ is a **Prime number**.
For $\boldsymbol{N} > 13$, $\boldsymbol{P = 6^{m+1} \bullet N - 1}$ is a **Prime number if and only if**
$\boldsymbol{N \neq i^{m+1}\, Mod(6i+1) + (6i+1)a}$ ; $i, a \in \boldsymbol{Z}^+$

## 3. Usage of $6^m$ Theorem

### 3.1 validating $6^m$ Theorem for m=1, i.e. $S1_{6^1}$

### 3.2 Mod theorem for $i^2/(Ai+B)$

When $i = Ap + q\ , p \in \{\mathbf{0}, \mathbf{Z}^+\}\ ,\ A \in \mathbf{Z}^+\ ,\ q = \{1,2,3,\ldots,A\}\ , Aq > B,$ then,

$$\boldsymbol{R = i^2 Mod(Ai+B) = iq - Bp\ ; p = \{0,1,2,3,\ldots\},\qquad q = \{1,2,3,\ldots,A\}}$$

From above equation, $i^2Mod(6i+1)$ given by,

$$i^2Mod(6i+1) = \{5p+1, 11p+4, 17p+9, 23p+16, 29p+25, 35p+36,, p \in \{\mathbf{0}, \mathbf{Z}^+\}\}$$

#### 3.2.1 Proving $\boldsymbol{R = i^2Mod(Ai+B) = iq - Bp\ ; p = \{0,1,2,3,\ldots\},\ q = \{1,2,3,\ldots,A\}}$

We can represent any number $i$ as,

$i = Ap + q\ ; p = \{0,1,2,3,\ldots\},\ A \in \mathbf{Z}^+,\ \ q = \{1,2,3,\ldots,A\}$

The quotient of $i^2/(Ai+B)$ is assumed as $p$ and the remainder as $R$ and then $i^2$ is,

$$i^2 = p(Ai+B) + R$$
$$(Ap+q)^2 = p[\,A(Ap+q)+B] + R$$
$$R = (Ap+q)^2 - p[\,A(Ap+q)+B]$$
$$R = (Ap+q)q - pB$$
$$R = iq - pB \ \ldots\ldots\ldots\ldots\ (7)$$
$$R = (Ap+q\,)q - pB$$
$$= (Aq-B)P + q^2$$
$$(Aq-B) > 0;\text{ then } R > 0$$

Then the quotient of $i^2/(Ai+B)$ is assumed as $(p+1)$ and the remainder as $R'$ and then $i^2$ is,

$$(Ap+q)^2 = (p+1)[\,A(Ap+q)+B] + R'$$
$$R' = (Ap+q)^2 - (p+1)[\,A(Ap+q)+B]$$
$$R' = (Ap+q)^2 - p[\,A(Ap+q)+B] - [\,A(Ap+q)+B]$$
$$R' = iq - pB - Ai - B$$
$$R' = i(q-A) - (p+1)B$$

$q \le A\ ;\ \ i,p,B > 0$ Then,

$R' < 0$ and therefore, the first assumption for the remainder $R = iq - pB$ is correct.

$R = i^2Mod(Ai+B) = iq - Bp\ ; p = \{0,1,2,3,\ldots\},\qquad q = \{1,2,3,\ldots,A\}$

### 3.3 Case study for $\boldsymbol{A = 6, B = 1}$

For $i^2Mod(6i+1)$ case;

$$A = 6, B = 1\ and, p = \{0,1,2,3,\ldots\},\qquad q = \{1,2,3,4,5,6\}$$

Therefore, $i^2Mod(6i+1) = iq - p\quad ; i = 6p+q$

$$i^2Mod(6i+1) = (6p+q)q - p$$

For

$$\left.\begin{array}{ll} q = 1, & R = i^2Mod(6i+1) = (6p+1)\bullet 1 - p = 5p+1 \\ q = 2, & R = i^2Mod(6i+1) = (6p+2)\bullet 2 - p = 11p+4 \\ q = 3, & R = i^2Mod(6i+1) = (6p+3)\bullet 3 - p = 17p+9 \\ q = 4, & R = i^2Mod(6i+1) = (6p+4)\bullet 4 - p = 23p+16 \\ q = 5, & R = i^2Mod(6i+1) = (6p+5)\bullet 5 - p = 29p+25 \\ q = 6, & R = i^2Mod(6i+1) = (6p+6)\bullet 6 - p = 35p+36 \end{array}\right\} \ldots\ldots\ldots\ldots\ (8)$$

Therefore, every single $i^2Mod(6i+1)$ number is given by the above set of equations in (8) .

## 3.4. Verifying $6^m$ Theorem for m=1

### 3.4.1 Verifying $6^m$ Theorem for m=1 and $\mathbf{1 < N \leq 13}$ $\boldsymbol{and}$ $\mathbf{N \neq 8}$

According to the solutions in ⑧, $i^2 Mod(6i+1) \neq 2 \neq N$ and as claimed by the **$6^m$ Theorem,** $P = 6^2 \bullet 2 - 1 = 71$ is a Prime number.
The same result can be shown for $i^2 Mod(6i+1) \neq 12 \neq N$ ; $P = 6^2 \bullet 12 - 1 = 431$ is also a Prime number.

### 3.4.2 Verifying $6^m$ Theorem for m=1 and $\boldsymbol{N} > 13$

From the solutions in ⑧, $i^2 Mod(6i+1) \neq 20 \neq N$ and then according to **$6^m$ Theorem,**
$\boldsymbol{N \neq i^{m+1}\ Mod(6i+1) + (6i+1)a}$ ; $i, a \in \boldsymbol{Z}^+$ must be satisfied for $\boldsymbol{N} > 13$.
Considering $N \neq i^2\ Mod(6i+1) + (6i+1)a$,
For $i = 1$; $i^2\ Mod(6i+1) = 1$ which leads to $20 = 1 + 7a$ , and $a \notin \boldsymbol{Z}^+$
Similarly,
For $i = 2$; $i^2\ Mod(6i+1) = 4$ , $20 = 4 + 13a$ , and $a \notin \boldsymbol{Z}^+$
For $i = 3$; $i^2\ Mod(6i+1) = 9$ and $20 = 9 + 19a$ , and $a \notin \boldsymbol{Z}^+$
For $i \geq 4$ ; $(6i+1)a > 20$ and therefore cannot find a condition satisfying $a \in \boldsymbol{Z}^+$
Hence, $i^2 Mod(6i+1) + (6i+1)a \neq 20 \neq N$ ; $P = 6^2 \bullet 20 - 1 = 719$ is a Prime number.